\documentclass[12pt]{article}

\usepackage{amsmath,amsfonts,amssymb,amsbsy}

\usepackage{color}

\newcounter{sect1}
\newcounter{subsect1}

\newfont{\gotikai}{eufm10}

\newcommand {\beq}{\begin{equation}}
\newcommand {\eeq}{\end{equation}}

\newtheorem{theorem}{\indent Theorem}
\newtheorem{corollary}{\indent Corollary}
\newtheorem{Lemma}{\indent Lemma}
\newtheorem{Remark}{\indent Remark}

\def\({\left(}
\def\){\right)}

\def\({\left(}
\def\){\right)}

\def\[{\left[}
\def\]{\right]}

\def\|{\left|}
\def\|{\right|}

\def\0{{\boldsymbol{0}}}

\def\00{{\boldsymbol{0}}}

\begin{document}

\begin{center}
\textbf{ Bounds for the Success Probability in the Odds Theorem}
\end{center}

\ \ \

\begin{center}
\textbf{A.M.~Kabaeva, A.V.~Logachov and A.A.~Yambartsev}
\end{center}

\ \ \

\textbf{Abstract}. Bruss's odds theorem \cite{Bruss1} addresses the problem of determining the optimal stopping time for sequences of independent indicator functions. In this note, we derive upper and lower bounds for the success probability under the optimal stopping rule. These bounds depend on the number of independent events under consideration and on a deterministic index specifying the stopping time. Moreover, the bounds are sharp: we provide explicit examples in which the corresponding inequalities are attained with equality.


\ \ \

\section{Introduction and main results}

In February 1960, Martin Gardner published the following problem (Problem 3) in his \textit{Mathematical Games} column in Scientific American, under the title \textit{A Fifth Collection of ``Brain-Teasers''} \cite{Gardner}. On page 150, he wrote:
\begin{enumerate}
    \item[] \textit{Ask someone to take as many slips of paper as he pleases, and on each slip write a different positive number. The numbers may range from small fractions of one to a number the size of a ``googol'' (1 followed by a hundred zeros) or even larger. These slips are turned face-down and shuffled over the top of a table. One at a time you turn the slips face-up. The aim is to stop turning when you come to the number that you guess to be the largest of the series. You cannot go back and pick a previously turned slip. If you turn over all the slips, then of course you must pick the last one turned.}
\end{enumerate}
In the English-speaking mathematical literature, this problem has become known as the famous \textit{secretary problem}, even though the word ``secretary'' does not appear in this Problem 3. 
After turning past two consecutive pages of advertisements (151 and 152), we reach page 153, where the text reads as follows:
\begin{enumerate}
    \item[]  \textit{The game has many interesting applications. For example, a bachelor girl decides to marry before the end of the year. She estimates that she will meet 10 men who can be persuaded to propose, but once she has rejected a proposal, the man will not try again. What strategy should she follow to maximize her chances of accepting the top man of the 10, and what is the probability that she will succeed?}
\end{enumerate}
This is likely the reason the problem is also referred to as the \textit{marriage problem}. It should be noted, however, that in contemporary English-language literature (including arXiv, as well as journals in algorithmic theory, operations research, and economics), the term ``marriage problem'' is almost exclusively used to denote the Stable Matching Problem. This problem concerns the construction of stable pairwise matchings between two sets (such as students and universities, or men and women), based on priorities and preferences. It is unrelated to the single-choice maximum selection problem commonly referred to as the secretary problem.

In the Soviet Union, and subsequently in the broader Russian-language popular science literature, this problem acquired a second designation, the ``marriage problem''. The literal rendering of its Russian counterpart is \textit{the problem of the choosy bride}, and this terminology has remained in use ever since.

For more historical issues about this problem, we refer the reader to \cite{Ferguson} and the good review \cite{Freem}. It is unclear whether a modern survey of this problem is feasible, given the vast and nearly unbounded literature on optimal stopping problems related to it and its generalizations. Such a survey would undoubtedly be of interest and is still awaiting its author.

The modern formulation of this problem no longer relies on uniformly random permutations; instead, it is posed (and solved) in a more general framework of a sequence of independent indicator variables (see \cite{Freem, Bruss1, Bruss3}), a formulation that we adopt in this paper. Consider mutually independent events $A_1, \dots, A_n$, $n \in \mathbb{N}$, defined on a probability space $(\Omega, \mathfrak{F}, \mathbf{P})$, and let their indicators be $I_j := \mathbf{I}(A_j)$ with success probabilities $p_j := \mathbf{P}(I_j=1)$, for $1 \leq j \leq n$. We assume that $n$ is fixed and the distributions of the random variables $I_j$, $1 \leq j \leq n$, are known. This corresponds to the so-called full-information framework. In the context of the secretary problem, $I_j=1$ means that the $j$th item is better than all of the previous $j-1$ items, and hence $p_j=\frac{1}{j}$.

The strategy is to first observe and reject the initial $s$ items (candidates), and then to accept the first one that is better than all previous items. We observe the indicator values sequentially, keeping track of all past and current values but without knowledge of future ones, and the goal is to stop at the last ``success'' (i.e., if $\tau$ is the stopping time, then the desired outcome is $I_\tau = 1, I_{\tau+1} = 0,\dots, I_n = 0$). If we succeed, we say that we ``win''. In the case of the secretary problem, Lindley \cite{Lindley} solved it by identifying the index $s$ as
$$
s= \arg \{k: a_{k-1} \ge 1 > a_k\},
$$
where
$$
a_k=\frac 1k + \frac{1}{k+1} + \dots + \frac{1}{n-1}.
$$

In the more general setting considered in this note, the solution was provided by Bruss \cite{Bruss1} in a paper whose title directly refers to the determination of the threshold step $s$: \textit{Sum the Odds to One and Stop}. Let $\mathfrak{L}$ denote the set of all stopping times, that is, a random variable $t \in \mathfrak{L}$ if
$\{ t = k \} \in \mathfrak{F}_k$ for every $1 \leq k \leq n$,
where $\mathfrak{F}_k$ is the $\sigma$-algebra generated by the random variables $I_1, \dots, I_k$. In \cite{Bruss1}, the problem of selecting an optimal stopping time was solved; namely, a random variable $\tau \in \mathfrak{L}$ such that
$$
\mathbf{P}\big(I_\tau = 1, \, I_{\tau+1} = 0, \dots, I_n = 0\big)
= \max_{t \in \mathfrak{L}} \,
\mathbf{P}\big(I_t = 1, \, I_{t+1} = 0, \dots, I_n = 0\big).
$$
For the reader's convenience, we now state this result. Let $I_1, \dots, I_n$ be independent indicator random variables, with $p_j = \mathbf{E} I_j$, $q_j = 1 - p_j$, and odds $r_j = \frac{p_j}{q_j}$, for $1 \leq j \leq n$.

\begin{theorem} \label{th1} (\cite{Bruss1}, \textbf{Bruss' odds theorem
of optimal stopping})
An optimal stopping rule $\tau \in \mathfrak{L}$ is to stop at the first index $k\ge s$ (if such exists) with $I_k = 1$, where
\begin{equation}\label{25-08-2025-1}
s = \max\left\{ 1,\  \max\left\{ l : \sum_{j=l}^n r_j \geq 1 \right\} \right\},
\end{equation}
with the convention $\max\{\varnothing\} = -\infty$. The corresponding probability of ``winning'' is
\begin{equation}\label{11.07.1}
V_n := \mathbf{P}(I_\tau = 1, I_{\tau+1} = 0, \dots, I_n = 0)
= \prod_{j=s}^n q_j \sum_{l=s}^n r_l.
\end{equation}
\end{theorem}



We are interested in upper and lower bounds for the success probability $V_n$. To begin with, we recall two well-known lower bounds from \cite{Bruss3} and \cite{AllIsl}.
\begin{enumerate}
    \item In \cite{Bruss3} Bruss showed that, if $\sum\limits_{j=1}^n r_j \geq 1$, then
    $$V_n > \dfrac{1}{e}.$$
    \item The work \cite{AllIsl} provides a generally sharp lower bound: if $\sum\limits_{j=1}^n r_j \geq 1$, then
    $$V_n \geq \left( 1 - \frac{1}{n+1} \right)^{n},$$
    with equality when $p_1 = \dots = p_n = \dfrac{1}{n+1}$.
\end{enumerate}




Unlike these known lower bounds, which depend only on $n$, our results concern both lower and upper bounds that also take into account the threshold index $s$ defined by \eqref{25-08-2025-1}, i.e., the index from which one stops at the first success. Each bound we provide in this paper is sharp: there are examples of sequences $I_j, \ 1\le j\le n$, for which the corresponding inequality is attained as an equality.


Before formulating the main results, we introduce the notation for the sum of odds:
$$
R_l:=\sum\limits_{j=l}^n r_j, \ \ \ 1\leq l \leq n.
$$
The main results of this paper are summarized in the following two theorems. Recall that the threshold index $s$ defined by \eqref{25-08-2025-1}, can now be expressed in the new notation as 
$$
s = \max\bigl\{ 1,\  \max\left\{ l : R_l \geq 1 \right\} \bigr\}.
$$

\begin{theorem}\label{th3}
\textbf{(Upper bound on the success probability in the odds theorem)}
The success probability \eqref{11.07.1} satisfies 
$$
V_n \leq \frac{R_s}{1 + R_s}.
$$
This upper bound is attained only for the sequence of indicator-function probabilities given by
\begin{equation}\label{11.07.8}
p_s = \frac{R_s}{1 + R_s}, \quad p_j = 0 \quad \text{for all } j \neq s.
\end{equation}
\end{theorem}

\begin{theorem}\label{th4}
\textbf{(Lower bounds on the success probability in the odds theorem)}
The following lower bounds hold depending on the value of $R_s$:

\begin{enumerate}
    \item If $0 \leq R_s < 1$ and $s=1$, then
    $$
    V_n \geq R_1 \left( 1 + \frac{R_1}{n} \right)^{-n}.
    $$
    Equality is attained when
    \begin{equation}\label{11.07.2}
    p_1 = \cdots = p_n = \frac{R_1}{n + R_1}.
    \end{equation}

    \item If $1 \leq R_s \leq 1 + \frac{1}{n - s}$, then
    $$
    V_n \geq \left( 1 + \frac{1}{n - s + 1} \right)^{-(n - s + 1)}.
    $$
    Equality holds if
    \begin{equation}\label{11.07.3}
    p_s = \cdots = p_n = \frac{1}{n - s + 2}.
    \end{equation}

    \item If $R_s > 1 + \frac{1}{n - s}$, then
    $$
    V_n > \left( 1 + \frac{1}{n - s} \right)^{-(n - s)}.
    $$
    By choosing the parameter $\alpha\in (0,1)$ sufficiently close to $1$ and setting
    \begin{equation}\label{11.07.5}
    p_s = \frac{1 + (2 - 2\alpha)(n - s)}{1 + (3 - 2\alpha)(n - s)}, \quad
    p_{s+1} = \cdots = p_n = \frac{\alpha}{n - s + \alpha},
    \end{equation}
    the success probability $V_n$ can be made arbitrarily close to $\left( 1 + \frac{1}{n - s} \right)^{-(n - s)}$.
\end{enumerate}
\end{theorem}

\begin{corollary}
If only $n$ and $s \geq 2$ are known and there is no information about the value of $R_s$, then the following lower bound holds:
$$
V_n \geq \left(1 + \frac{1}{n - s + 1}\right)^{-(n - s + 1)}.
$$
\end{corollary}

\begin{Remark}
Theorem~\ref{th3} and Theorem~\ref{th4} imply that, if $R_s=1$, then
$$
\frac{1}{e}<
\left(1+\dfrac{1}{n}\right)^{-n}\leq\left(1+\dfrac{1}{n-s+1}\right)^{-(n-s+1)}\!\!\!\! \leq V_n\leq\frac{1}{2}.
$$
\end{Remark}

\begin{Remark}
The inequalities in cases 2 and 3 depend on the value of $s$. As $s$ increases, the lower bound for the success probability under the optimal stopping strategy increases. In particular, setting $s=1$ in the lower bound of case 2 yields the lower bound from \cite{AllIsl} mentioned before.
\end{Remark}

\begin{Remark}
The main difference between cases 2 and 3 is that, in case 3, we always have $s < n$ and the equality
$$
r_s = \cdots = r_n = \frac{R_s}{n - s + 1}
$$
is impossible. Indeed, if this equality were satisfied, one could easily check that $R_{s+1}>1$ in this case:
$$
R_{s+1} = (n-s) \frac{R_s}{n - s + 1} > (n-s) \frac{1+\frac{1}{n-s}}{n - s + 1} = 1,
$$
which contradicts the definition of $s$.

Consequently, the proof of the lower bound in case 3 is reduced to finding of infimum of a function of two variables, $R_s$ and $R_{s+1}$, on the domain
$$
R_s > 1 + \frac{1}{n - s}, \quad 0 \leq R_{s+1} < 1,
$$
whereas the lower bound in case 2 is obtained by finding the infimum of a function of one variable, $R_s$, on the interval
$$
1 \leq R_s \leq 1 + \frac{1}{n - s}.
$$
\end{Remark}


\section{Proof of main results}

In all proofs, we regard the success probability $V_n$ as the function $V(r_s, \ldots, r_n)$ of $n - s + 1$ variables on the domain $r_j \geq 0$, $s \leq j \leq n$, subject to the constraint $\sum_{j=s}^n r_j = R_s$. Standard methods of calculus are then used to determine the maximum and minimum of $V$ on this domain.

\subsection{Proof of Theorem~\ref{th3}}

The success probability $V_n$, see \eqref{11.07.1}, can be represented in the following alternative form  
\begin{equation}\label{11.07.7}
V_n \equiv V(r_s, \ldots, r_n) := \frac{R_s}{\prod_{j=s}^n (1 + r_j)}.
\end{equation}
Using the exp–log representation, the function can be written as
$$
V(r_s, \ldots, r_n) = R_s \exp\left\{ - \sum_{j=s}^n \ln(1 + r_j) \right\},
$$
The upper bound \eqref{08.07.1} from Lemma~\ref{l1} (see below) implies the upper bound for the success probability
$$
V(r_s, \ldots, r_n) \leq R_s \exp\left\{ - \ln\left(1 + \sum_{j=s}^n r_j \right) \right\} = \frac{R_s}{1 + R_s}.
$$
Note that the equality holds if and only if 
$$
\prod_{j=s}^n (1 + r_j) = 1 + R_s.
$$
It is straightforward to observe that this occurs only when $r_s = R_s$ and $r_j = 0$ for all $j \neq s$, which corresponds precisely to (\ref{11.07.8}). This finishes the proof of the theorem. $\Box$

\vspace{0.5cm}

To complete the proof, we now establish Lemma~\ref{l1}.

\begin{Lemma}\label{l1}
For fixed two numbers $k\in \mathbb{N}$ and $a\geq 0$, consider the following domain
$$
D_a := \left\{(x_1,\dots,x_k)\in \mathbb{R}^k : x_j \ge 0,~1\le j \le k,~\sum\limits_{j=1}^k x_j = a \right\}.
$$
Then the following upper bound holds,
\begin{equation}\label{08.07.1}
\sup\limits_{(x_1,\dots,x_k)\in D_a} \left(-\sum\limits_{j=1}^k \ln(1+x_j)\right) \leq -\ln(1+a).
\end{equation}
Equality in (\ref{08.07.1}) is attained when $x_i = a$ and $x_j = 0$ for all $j \neq i$.
\end{Lemma}

\noindent
\textbf{Proof of Lemma~\ref{l1}:} The final claim of the lemma can be verified by noting that if $x_i=a$ for some $i$, and $x_j=0$ for all $j\neq i$, then 
$$
-\sum\limits_{j=1}^k\ln(1+x_j)=-\ln(1+a).
$$
We will prove (\ref{08.07.1}) by induction on $k, k\ge 1$. For the base case $k = 1$, the statement is immediate: the domain is $D_a:=\left\{x_1\in \mathbb{R}:x_1=a\right\}$ and for any $a\geq0$ we have
$$\sup\limits_{x_1\in D_a}\left(-\ln(1+x_1)\right)=-\ln(1+a).$$
Assume now, that (\ref{08.07.1}) holds for $k=l$, i.e., for any $a\geq 0$ the inequality 
$$
\sup\limits_{(x_1,\dots,x_l)\in D_a}\left(-\sum\limits_{j=1}^l\ln(1+x_j)\right)\leq-\ln(1+a)$$ is satisfied.
We will show that the same inequality (\ref{08.07.1}) holds for $k = l+1$ and any $a \ge 0$.

Using the inductive assumption, for any $a \ge 0$ and $(x_1,\dots,x_{l+1}) \in D_a$ we have
\beq\label{08.07.2}
\begin{aligned}
-\sum\limits_{j=1}^{l+1}\ln(1+x_j) &=-\sum\limits_{j=1}^l\ln(1+x_j)-\ln(1+x_{l+1}) \\
&\leq -\ln\left(1+\sum\limits_{j=1}^lx_j\right)-\ln(1+x_{l+1}).
\end{aligned}
\eeq

Denote $b:=\sum_{j=1}^lx_j$. The supremum on the right-hand side of inequality \eqref{08.07.2} is attained at the largest value of the function
$$f(b):=-\ln\left(1+b\right)-\ln(1+a-b)$$
on the interval $b\in[0,a]$. Elementary calculus shows that $f$ attains its maximum over $[0,a]$ at the endpoints,
\beq\label{08.07.3}
\sup_{b\in [0,a]} f(b) \le f(0) = f(a) = - \ln(1+a).
\eeq

Thus, (\ref{08.07.2}) together with \eqref{08.07.3} yields 
$$
\sup \left(-\sum_{j=1}^{l+1} \ln(1 + x_j)\right) \leqslant -\ln(1 + a),
$$
for any $a\geq0$, $(x_1,\dots,x_{l+1})\in D_a$.
This completes the proof of the lemma. $\Box$

\subsection{Proof of Theorem~\ref{th4}}

\textbf{Proof of item 1.} Consider the case $0\leqslant R_s<1$. In this case, by definition of index $s$, this implies that $s=1$. Using the fact that the geometric mean is less than the arithmetic mean, we obtain
$$
V_n=\frac{R_1}{\prod\limits_{j=1}^n(1+r_j)}
\geq R_1\left(\frac{1}{n}\sum\limits_{j=1}^n(1+r_j)\right)^{-n}
= R_1\left(1+\frac{R_1}{n}\right)^{-n}.
$$
It remains to note that equality is attained if
$$
r_1=\ldots=r_n=\frac{R_1}{n},
$$
which is equivalent to the condition (\ref{11.07.2}). $\Box$

\vspace{0.5cm}
\noindent
\textbf{Proof of item 2.}
Now, consider the case $1\leq R_s\leq 1+\frac{1}{n-s}$. Using again the inequality between the geometric and arithmetic means we obtain
\begin{equation}\label{th02}
V_n=\frac{R_s}{\prod\limits_{j=s}^n(1+r_j)}\geq R_s\left(1+\frac{R_s}{n-s+1}\right)^{-(n-s+1)}.
\end{equation}
Taking the logarithm of the right-hand side of inequality \eqref{th02}, we obtain
$$
\begin{aligned}
g(R_s):=&\ln\left(R_s\left(1+\frac{R_s}{n-s+1}\right)^{-(n-s+1)}\right) \\
=&\ln R_s-(n-s+1)\ln \left(1+\frac{R_s}{n-s+1}\right).
\end{aligned}
$$
Let us find the point at which the function $g$ attains its minimum value on the interval  $1\leq R_s\leq 1+\frac{1}{n-s}$. Observe that the following inequality holds on this interval
$$
\begin{aligned}
g'(R_s)&=\frac{1}{R_s}-\frac{n-s+1}{n-s+1+R_s} \\
&=\frac{n-s+1+R_s-(n-s+1)R_s}{R_s(n-s+1+R_s)}
=\frac{n-s+1-(n-s)R_s}{R_s(n-s+1+R_s)}\geq 0.
\end{aligned}
$$
Thus, the function $g$ monotonically increases with $R_s$ over the considered interval, and its minimum is attained at $R_s = 1$. Substituting this value into \eqref{th02}, we obtain
$$
V_n\geq\left(1+\frac{1}{n-s+1}\right)^{-(n-s+1)}.
$$
It is easy to see that equality is attained when
$$
r_s=\dots=r_n=\frac{1}{n-s+1},
$$
which implies condition (\ref{11.07.3}):
$$
p_i=\frac{r_i}{1+r_i} = \frac{1}{n-s+2}, \ \ \ i=s, \ldots, n.
$$
$\Box$

We close this section with the following remark. The distinction in the values of $R_s$ between cases 2 and 3 in the theorem arises from the equation $g'(R_s)=0$. This equation is equivalent to $$n-s+1-(n-s)R_s=0,$$ which yields the separating value $$R_s=1+\frac{1}{n-s}.$$ 

\vspace{0.5cm}
\noindent
\textbf{Proof of item 3.} Consider now the case $R_s> 1+\frac{1}{n-s}$. Once again, using the fact that the geometric mean does not exceed the arithmetic mean, we obtain
\begin{equation}\label{t21}
\begin{aligned}
V_n =\frac{R_s}{\prod\limits_{j=s}^n(1+r_j)}&=\frac{R_s}{(1+R_s-R_{s+1})\prod\limits_{j=s+1}^n(1+r_j)} \\
&\geq\frac{R_s}{(1+R_s-R_{s+1})}\left(\frac{1}{n-s}\sum\limits_{j=s+1}^n(1+r_j)\right)^{-(n-s)}
\\
&=\frac{R_s}{(1+R_s-R_{s+1})}\left(1+\frac{R_{s+1}}{n-s}\right)^{-(n-s)}.
\end{aligned}
\end{equation}
Taking the logarithm of the final expression for the lower bound on the right-hand side of inequality (\ref{t21}), we obtain
$$
\begin{aligned}
f(R_s, R_{s+1}):=&\ln\Bigg(\frac{R_s}{1+R_s-R_{s-1}}\bigg(1+\frac{R_{s+1}}{n-s}\bigg)^{-(n-s)}\Bigg) \\
=&\ln R_s-\ln(1+R_s-R_{s+1})-(n-s)\ln\left(1+\frac{R_{s+1}}{n-s}\right).
\end{aligned}
$$
We find the infimum of the function $f$ in the domain
$$
D:=\left\{(R_s,R_{s+1})\in \mathbb{R}^2:~ R_s> 1+\frac{1}{n-s},  \  0\leq R_{s+1}<1\right\}.
$$
Let us find the partial derivatives of the function $f$
$$
\frac{\partial{f}}{\partial{R_s}}=\frac{1}{R_s}-\frac{1}{1+R_s-R_{s+1}}, \ \ \  \frac{\partial{f}}{\partial{R_{s+1}}}=\frac{1}{1+R_s-R_{s+1}}-\frac{n-s}{n-s+R_{s+1}}.
$$
It is easy to see that $\frac{\partial{f}}{\partial{R_s}}>0$ on the domain $D$. Therefore, $f$ is a monotonically increasing function of $R_s$ on $D$.

Let us bound $\frac{\partial f}{\partial R_{s+1}}$ from above. On the given domain $D$, we have the following
$$
\begin{aligned}
\frac{\partial f}{\partial R_{s+1}} &=\frac{n-s+R_{s+1}-(n-s)(1+R_s-R_{s+1})}{(1+R_s-R_{s+1})(n-s+R_{s+1})} \\
& <\frac{n-s+R_{s+1}-(n-s)R_s}{(1+R_s-R_{s+1})(n-s+R_{s+1})}<\frac{n-s+R_{s+1}-(n-s+1)}{(1+R_s-R_{s+1})(n-s+R_{s+1})} \\
&=\frac{R_{s+1}-1}{(1+R_s-R_{s+1})(n-s+R_{s+1})}<0.
\end{aligned}
$$
Thus, the function $f$ decreases monotonically in $R_{s+1}$ on the domain $D$. Therefore, the infimum of the function $f$ on the considered domain is attained at
$$
R_s\downarrow 1+\frac{1}{n-s}, \ R_{s+1}\uparrow1.
$$

Substituting the obtained preliminary values into the right-hand side of inequality \eqref{t21}, we obtain
$$
V_n>\left(1+\frac{1}{n-s}\right)^{-(n-s)}.
$$
Choosing the parameter $\alpha\in (0,1)$ sufficiently close to $1$, and setting
$$
r_s=2-2\alpha+\frac{1}{n-s}, \ r_{s+1}=\dots=r_n=\frac{\alpha}{n-s},
$$
(which satisfies the relations \eqref{11.07.5}),
we will obtain values of $V_n$ arbitrarily close to $$\left(1+\frac{1}{n-s}\right)^{-(n-s)}.$$ This concludes the proof of the third item, and hence the proof of the theorem. $\Box$

\subsection*{Acknowledgments}

This work was supported by the Russian Science Foundation (RSF), grant no. 24-28-01047. 

\end{document}